\documentclass[12pt]{amsart}
\usepackage{amsthm}
\usepackage{amsmath, mathrsfs}
\usepackage{amssymb}
\usepackage{hyperref}
\usepackage{enumerate}
\usepackage{latexsym,array}
\usepackage{amsfonts}
\usepackage{shadow}

\newcommand{\bigA}{\mathcal{A}}

\newcommand{\bigS}{\mathcal{S}}

\newcommand{\bigP}{\mathcal{P}}
\newcommand{\bigM}{\mathcal{M}}

\newtheorem{Pa}{Paper}[section]
\newtheorem{Tm}[Pa]{{\bf Theorem}}
\newtheorem{La}[Pa]{{\bf Lemma}}
\newtheorem{Dn}[Pa]{{\bf Definition}}
\newtheorem{Cy}[Pa]{{\bf Corollary}}

\newtheorem{Rk}[Pa]{{\bf Remark}}
\newtheorem{Pn}[Pa]{{\bf Proposition}}

\newtheorem{Ex}[Pa]{{\bf Example}}

\date{}

\thanks{The authors thank the referee for her/his remarks and in particular
for the suggestion to consider the non-unimodular case. D. Alpay
thanks the Earl Katz family for endowing the chair which
supported his research, and the Binational Science Foundation
Grant number 2010117.
}
\author{Daniel Alpay}
\address{(DA) Department of Mathematics
\newline
Ben Gurion University of the Negev \newline P.O.B. 653,
\newline
Be'er Sheva 84105, \newline ISRAEL}
\email{dany@math.bgu.ac.il}
\author{Guy Salomon}
\address{(GS) Department of Mathematics
\newline
Ben Gurion University of the Negev \newline P.O.B. 653,
\newline
Be'er Sheva 84105, \newline ISRAEL} \email{guysal@math.bgu.ac.il}
\title[Topological convolution algebras]
{Topological convolution algebras}
\begin{document}
\maketitle


\begin{abstract}
In this paper we introduce a dual of reflexive Fr\'echet counterpart of
Banach algebras of the form $\bigcup_{p\in\mathbb N}\Phi_p^\prime$
(where the $\Phi_p^\prime$ are (dual of)  Banach spaces
with associated norms $\|\cdot\|_p$), which carry inequalities of the form
$\|ab\|_{p} \leq A_{p,q}\|a\|_{q}\|b\|_{p}$ and $\|ba\|_{p}
\leq A_{p,q}\|a\|_{q}\|b\|_{p}$ for $p> q+d$, where $d$ is preassigned and
$A_{p,q}$ is a constant.
We study the functional calculus and the spectrum of the elements
of these algebras.
We then focus on the particular case
$\Phi_p^\prime=L_2(S,\mu_p)$, where $S$ is a Borel semi-group in a locally
compact group $G$, and multiplication is convolution.
We give a sufficient
condition on the measures $\mu_p$ for such inequalities to hold.
Finally we
present three examples, one is the algebra of germs of holomorphic functions in zero,
the second related to Dirichlet series and the third in the setting of
non commutative stochastic distributions.
\end{abstract}

\section{Introduction}
\setcounter{equation}{0}
Let $G$ be a locally compact topological group with a left
Haar measure $\mu$. The convolution of two measurable
functions $f$ and $g$ is defined by
\[
(f*g)(x)=\int_{G}f(y)g(y^{-1}x)d\mu(y).
\]
It is well known that $L_1(G,\mu)$ is a Banach algebra with the
convolution product, while $L_2(G,\mu)$
is usually not closed under the convolution. More precisely, Rickert
showed in 1968 that $L_2(G,\mu)$ is closed under convolution if and
only if $G$ is compact; see \cite{MR0228930}.
In this case $G$ is unimodular (i.e. its left and right Haar measure
coincide) and it holds that
\[
\|f*g\| \leq \sqrt{\mu(G)}\|f\|\|g\|,\quad \text{ for any }f,g \in L_2(G,\mu).
\]
Thus, $L_2(G,\mu)$ is a convolution algebra which is also a Hilbert space.
In this paper we introduce and study a new type of
convolution algebras which behave locally as Hilbert spaces, rather
than being Banach spaces, even when the group $G$ is
not compact. More precisely,
let $(\mu_p)$ be a sequence of measures on $G$ such that
\[
\mu \gg \mu_1 \gg \mu_2 \gg \cdots,
\]
(where $\mu$ is the left Haar measure of $G$) and let $S \subseteq G$ be a
Borel semi-group with the following property: There exists a non negative integer $d$
such that for any $p> q+d$
and any $f \in L_2(S,\mu_q)$ and $g \in L_2(S,\mu_p)$, the products $f*g$
and $g*f$ belong to $L_2(S,\mu_p)$ and
\begin{equation}
\label{ineq}
\|f*g\|_p \leq A_{p,q}\|f\|_q\|g\|_p \quad \text{ and } \quad \|g*f\|_p \leq A_{p,q}\|f\|_q\|g\|_p,
\end{equation}
where $\|\cdot\|_p$ is the norm associate to $L_2(S,\mu_p)$, where $A_{p,q}$
is a positive
constant and where the convolution of two measurable functions with respect to
$S$ is defined by
\eqref{convS}.\\
%

More generally, we consider dual of reflexive Fr\'echet spaces of
the form $\bigcup_{p\in\mathbb N}\Phi_p^\prime$, where the
$\Phi_p^\prime$ are (dual of) Banach spaces with associated norms
$\|\cdot\|_p$), which are moreover algebras and carry
inequalities of the form
\[
\|ab\|_{p} \leq A_{p,q}\|a\|_{q}\|b\|_{p} \quad \text{ and }\quad \|ba\|_{p} \leq A_{p,q}\|a\|_{q}\|b\|_{p}
\]
for $p> q+d$, where $d$ is preassigned and $A_{p,q}$ is a constant. We
call these spaces {\sl strong algebras}.
They are topological algebras (see Theorem \ref{SAcont}). Furthermore,
the multiplication operators
\[
L_a \,:\, x \mapsto ax \text{ and } R_a\,:\, x \mapsto xa
\]
are bounded from the Banach space $\Phi_p^\prime$ into
itself where $a\in \Phi_q^\prime$ and $p>q+d$.
Such a setting allows to consider power series. If $\sum_{n=0}^\infty c_n z^n$
converges in the open disk with
radius $R$, then for any $a \in \Phi_q^\prime$ with $\| a\|_{q}<\frac{R}{A_{q+d+1,q}}$,  we obtain
\[
\sum_{n=0}^\infty |c_n|\| a^{n}\|_{q+d+1} \leq \sum_{n=0}^\infty |c_n|
(A_{p,q}\| f\|_{q})^n<\infty,
\]
and hence $\sum_{n=0}^\infty c_n a^{n} \in \Phi^\prime_{q+d+1}$. In
this way we are also able to consider the invertible
elements of the algebra $\bigcup_p \Phi_p^\prime$.\\

When we return to the case where $\Phi_p^\prime=L_2(S,\mu_p)$,
we give a sufficient condition on the sequence of measures $(\mu_p)$ such that
\eqref{ineq} holds.
More precisely, we show that if their Radon-Nikodym derivatives with respect
to the left Haar measure are submultiplicative, and if there exists $d$
such that for any $p>q+d$, the functions $\frac{d\mu_p}{d\mu_q}$ belong to
$L_1$ of both the left and right Haar measures, then \eqref{ineq} holds
(see Theorem \ref{DFsemi}).\\

There is one well known example for such an algebra, namely the germs of
holomorphic functions at the origin, with the pointwise multiplication
(which is in terms of power series, a convolution).
See e.g. \cite{MR745622}. We show that it can be identified as a union of Hardy
spaces on decreasing sequence of disks, and that with respect to the associated
norms it satisfies \eqref{ineq} (see Theorem \ref{SAgerms}). This
allows to develop analytic calculus of germs, and to conclude easily
some topological properties of this space, e.g. that it is nuclear
(see Corollary \ref{TopGerms}).\\

Another example for a convolution algebra associated to a
semi-group and which satisfies \eqref{ineq} is the space $\bigP$
of all functions $f:[1,\infty) \to \mathbb R$ such that there
exists $p \in \mathbb N$ such that
\[
\int_1^\infty |f(x)|^2 \frac{dx}{x^{p+1}}<\infty,
\]
with the convolution
\[
(f*g)(x)=\int_1^xf(y)g\left(\frac{x}{y}\right)\frac{dy}{y}, \quad \forall x \in [1,\infty).
\]
Inequalities \eqref{ineq} lead then to apparently new inequalities associated to
Dirichlet series.\\

Yet another example (which was the motivation for the present work) is given by
the Kondratiev space of Gaussian stochastic distributions. It can be
defined as $\bigcup_{p\in\mathbb N} L_2(\ell,\mu_p)$ where $\ell$
is the free commutative semi-group generated by the natural
numbers, and
\[
\mu_p(\alpha)=(2\mathbb N)^{-\alpha p}=\prod_{n =1}^\infty (2n)^{-\alpha(n)\cdot p}, \quad \forall \alpha \in \ell.
\]
In this case, the convolution is sometimes called Wick product. In 1996, V\aa ge (see \cite{vage96},
\cite[p. 118]{MR2387368}) showed that the
Kondratiev space with the convolution product satisfies \eqref{ineq}, where
\[
A_{p,q}=\left(\sum _{\alpha \in \ell} (2 \mathbb{N})^{-(p-q)
\alpha}\right)^{\frac12}<\infty.
\]
This fact plays a key role in stochastic partial differential
equations and in stochastic linear systems theory; see
\cite{MR1408433,alp,MR2610579}.
We show here that the non-commutative version of this space still satisfies
\eqref{ineq}.\\

The outline of the paper is as follows: In Section 2 we study
topological and spectral properties of strong algebras and of
their elements. In Section 3 we consider the case where the
multiplication is a convolution. The above mentioned examples are
presented in Sections 4, 5 and 6.



\section{Strong algebras}
\setcounter{equation}{0}

We now introduce a family of locally convex algebras, more precisely duals
of reflexive complete countably normed spaces,
which satisfy special inequalities. We recall that a countably normed
space $\Phi$ is a locally convex space whose topology is defined using a
countable set of compatible norms $(\|\cdot\|_n)$  i.e. norms such that if a
sequence $(x_n)$ that is a Cauchy sequence in the norms $\|\cdot\|_p$ and
$\|\cdot\|_q$ converges to zero in one of these norms, then it also converges
to zero in the other. The sequence of norms $(\|\cdot\|_n)$ can be always assumed
to be non-decreasing.
Denoting by $\Phi_p$ the completion of $\Phi$ with respect to the norm $\|\cdot\|_p$,
we obtain a sequence of Banach spaces
\[
\Phi_1 \supseteq \Phi_2 \supseteq \cdots \supseteq \Phi_p \supseteq \cdots.
\]
It is a well known result that $\Phi$  is complete if and only if $\Phi=\bigcap
\Phi_p$, and $\Phi$ is a Banach space if and only if there exists some $p_0$ such
that $\Phi_p=\Phi_{p_0}$ for all $p>p_0$.
Denoting by $\Phi'$ the dual space of $\Phi$ and by $\Phi_p'$ the dual of $\Phi_p$
it is clear that
\[
\Phi_1' \subseteq \Phi_2' \subseteq \cdots \subseteq \Phi_p' \subseteq \cdots,
\]
and that $\Phi'=\bigcup \Phi_p'$. A functional in $f\in\Phi_p'$ has the respective norm $\|f\|_p=\sup_{\|x\|_p\leq 1}|f(x)|$, and these norms on $\Phi'$ form a decreasing sequence. For further reading on countably
normed spaces and their duals we refer to \cite{GS2_english}.

\begin{Dn}
Let $\bigA=\bigcup_p \Phi_p'$ be a dual of a complete reflexive countably normed space. We call $\bigA$ a strong algebra if
it satisfies the property that there exists a constant $d$ such that for any $q$ and for any $p>q+d$ there exists a positive constant $A_{p,q}$ such that for any $a \in \Phi_q'$ and $b \in \Phi_p'$,
\[
\|ab\|_{p} \leq A_{p,q}\|a\|_{q}\|b\|_{p} \quad \text{ and } \quad \|ba\|_{p} \leq A_{p,q}\|a\|_{q}\|b\|_{p} .
\]
\end{Dn}

We topologized $\bigA$ with the strong topology, that is, a neighborhood of zero is defined by means of any bounded set $B \subseteq \bigcap_p \Phi_p$ and any number $\epsilon>0$, as the set of all $a \in \bigA$ for which
\[
\sup_{b\in B}|a(b)|<\epsilon.
\]
With this topology $\bigA$ is locally convex. Since $\bigA$ is the dual of
the reflexive Fr\'echet space $\bigcap \Phi_p$, its topology coincides with its topology as the inductive limit of the Banach spaces $\Phi_p'$ (see the proof of \cite[IV.23, Proposition 4]{BourbakiTVS}). In particular, it satisfies the universal property of inductive limits, i.e. any linear map from $\bigA$ to another locally convex space is continuous if and only if the restriction of the map to the any of the
spaces $\Phi_p'$  is continuous (see \cite[II.29]{BourbakiTVS}).\\

We now show that the multiplication is continuous in $\bigA$.
Before that, we show it is separately continuous.

\begin{Pn}
\label{pwcont}
Let $a \in \bigA$.
Then the linear mappings $L_a:x \mapsto ax$, $R_a:x \mapsto xa$ are continuous.
\end{Pn}
\begin{proof}[Proof]
Suppose that $a \in \Phi_q'$, and
let $L_a|_{\Phi_r'} : \Phi_r' \to  \bigA$ be the restriction of the
map $L_a$ to $\Phi_r'$.
If $B$ is a bounded set of $\Phi_r'$ then in particular we may choose
$p \geq q+d$ such that $p \geq r$, so
$B \subseteq \{ x \in \Phi_p': \|x\|_p <\lambda \}$.
Thus, for any $x \in B$
\[
\|L_a|_{\Phi_p'} (x) \|_q \leq A_{p,q} \lambda \|x\|_p.
\]
Hence, $L_a|_{\Phi_p'}(B)$ is bounded in $\Phi_p'$ and hence in $\bigA$.
Thus, for any $r$, $L_a|_{\Phi_r'} :\Phi_r' \to \bigA$ is bounded and hence continuous.
Since $\bigA=\bigcup_{p \in \mathbb N} \Phi_p'$ is a strong dual of the
reflexive Fr\'echet space $\bigcap \Phi_p$, it is the inductive limit of
the Hilbert spaces $\Phi_p'$.
So by the universal property of inductive limits, $L_a$ is continuous.
The proof for $R_a$ is similar.
\end{proof}

\begin{Tm}\label{SAcont}
Let $\bigA$ be a strong algebra. Then the multiplication is a continuous
function $\bigA
\times \bigA\to \bigA$ in the
strong topology. Hence $(\bigA,+,\cdot)$ is a
topological $\mathbb C$-algebra.
\end{Tm}

This follows immediately from Proposition \ref{pwcont} together with the
following theorem, proved in
\cite[IV.26]{BourbakiTVS}.
\begin{Tm}
Let $E_1$ and $E_2$ be two reflexive Fr\'echet spaces, and let $G$ a
locally convex Hausdorff space. For $i=1,2$, let $F_i$ be the strong dual
of $E_i$. Then every separately continuous bilinear mapping $u:F_1 \times F_2
\to G$ is continuous.
\end{Tm}

Henceforward, we assume that $\bigA$ is a unital strong algebra.
The following theorems shows that, as in the Banach algebra case,
one can develop an analytic calculus in strong algebras.

\begin{Tm}\label{power}
Assuming $\sum_{n=0}^\infty c_n z^n$ converges in the open
disk with radius $R$, then for any $a \in \bigA$ such that there
exist $p,q$ with $p>q+d$ and $A_{p,q}\| a\|_{q}<R$ it holds that
\[
\sum_{n=0}^\infty c_n a^{n} \in \Phi_p' \subseteq \bigA.
\]
\end{Tm}
\begin{proof}[Proof]
This follows from
\[
\sum_{n=0}^\infty |c_n|\| a^{n}\|_{p} \leq \sum_{n=0}^\infty |c_n|
(A_{p,q}\| a\|_{q})^n<\infty.
\]
\end{proof}

\begin{Tm}\label{invert}
Let $a \in \bigA$ be such that there exist $p,q$ such that $p>q+d$ and $A_{p,q}\|a\|_q<1$ then $1-a$ is invertible
(from both sides) and it holds that
\[
\|(1-a)^{-1}\|_{p} \leq \frac{1}{1-A_{p,q}\|a\|_q}, \quad \|1-(1-a)^{-1}\|_{p} \leq \frac{A_{p,q}\|a\|_q}{1-A_{p,q}\|a\|_q}.
\]
\end{Tm}
\begin{proof}[Proof]
Due to Theorem \ref{power} we have that
\[
\sum_{n=0}^\infty a^{n} \in \Phi_p' \subseteq \bigA.
\]
Moreover, clearly
\[
(1-a)\left(\sum_{n=0}^\infty a^{n}\right)=\left(\sum_{n=0}^\infty a^{n}\right)(1-a)=1,
\]
and we have that
\[
\|(1-a)^{-1}\|_p \leq \sum_{n=0}^\infty \| a^{n}\|_{p} \leq \sum_{n=0}^\infty  (A_{p,q}\| a\|_{q})^n=\frac{1}{1-A_{p,q}\|a\|_q},
\]
and
\[
\|1-(1-a)^{-1}\|_p \leq \sum_{n=1}^\infty \| a^{n}\|_{p} \leq \sum_{n=1}^\infty  (A_{p,q}\| a\|_{q})^n
=\frac{A_{p,q}\|a\|_q}{1-A_{p,q}\|a\|_q}.
\]
\end{proof}

As was mentioned before, since $\bigA=\bigcup_p \Phi_p'$ is the
strong dual of a reflexive Fr\'echet space, it is also the
inductive limit of the Banach spaces $\Phi_p'$, i.e. its strong
topology coincides is the finest locally-convex topology such
that the embeddings $\Phi_p' \hookrightarrow \bigA$ are
continuous. The following theorem refers to the case where the
strong topology of $\bigA$ is the finest topology (rather than
the finest {\em locally-convex} topology) such that these
mappings are continuous. There are two important cases when this
happens: the first is when $\bigA$ is a Banach algebra, and the
second is when the embeddings $\Phi_p' \hookrightarrow
\Phi_{p+1}'$ are compact. In particular, the examples described
in Sections 4 and 6 pertain to the second case (see Theorem
\ref{nuclear}).
\begin{Tm}\label{cycontinuous}
If the topology of $\bigA$ is the finest topology such that
the embeddings $\Phi_p' \hookrightarrow \bigA$ are continuous, then the set
of invertible elements $GL(\bigA)$ is open, and the function
$a \mapsto a^{-1}$ is continuous.
\end{Tm}

\begin{proof}[Proof]
Note that in this case, $U$ is open in $\bigA$ if and only if $U \cap \Phi_p'$ is open in $\Phi_p'$ for every $p$.
Let $U$ be an open set of $\bigA$, and let $p \in \mathbb N$.
Let $U_a$ be the set of all $b \in \bigA$ such that there exists $p>q+d$ for which
\[
\|b\|_{p} < \frac{1}{A_{p+d,p}A_{p,q}\|a^{-1}\|_q}.
\]
Clearly $U_a \cap \Phi_p'$ is open in $\Phi_p'$ for any $p$, so $U_a$ is open. Moreover, for any $b\in U_a$
\[
A_{p+d,p}\|a^{-1}b\|_{p} \leq A_{p,q}\|a^{-1}\|_q \|b\|_p<1.
\]
In view of Theorem \ref{invert}, $1-a^{-1}b$ is invertible, and therefore $a-b=a(1-a^{-1}b)$ is invertible too. Thus, $a+U_a \subseteq GL(\bigA)$, and so $GL(\bigA)$ is open.\\
Now, we note that,
\[
\begin{split}
(a+b)^{-1} -a^{-1}
&=\left(a(1+a^{-1}b)^{-1}\right)-a^{-1}\\
&=(1+a^{-1}b)^{-1}a^{-1}-a^{-1}\\
&=\left((1+a^{-1}b)^{-1}-1\right)a^{-1}.
\end{split}
\]
Therefore, for any $b \in U_a$,
\[
\begin{split}
\|(a+b)^{-1} -a^{-1}\|_{p+d}
&\leq A_{p+d,q}\|a^{-1}\|_q\|(1+a^{-1}b)^{-1}-1\|_{p+d}\\
&\leq A_{p+d,q}\|a^{-1}\|_q\frac{A_{p+d,p}\|a^{-1}b\|_{p}}{1-A_{p+d,p}\|a^{-1}b\|_p}\\
&\leq A_{p+d,q}\|a^{-1}\|_q\frac{A_{p+d,p}A_{p,q}\|a^{-1}\|_{q}\|b\|_{p}}{1-A_{p+d,p}A_{p,q}\|a^{-1}\|_q\|b\|_p}
\end{split}
\]
Thus, the function
\[
u:b \mapsto (a+b)^{-1} -a^{-1}
\]
satisfy $u(U_a \cap \Phi_p') \subseteq \Phi_{p+d}'$, and $u|_{U_a \cap \Phi_p'}$ is continuous with respect to the topologies of $U_a \cap \Phi_p'$ in the domain and $\Phi_{p+d}'$ in the range. \\
Now, let $V$ be an open set of $\bigA$.
Then, $u^{-1}(V) \cap \Phi_p' =u|_{U_a\cap \Phi_p'}^{-1}(V \cap \Phi_{p+d}')$ is open in $U_a \cap\Phi_p'$.
In particular, $u^{-1}(V) \cap \Phi_p'$ is open in $\Phi_p'$.
Since $p$ was arbitrary, $u^{-1}(V)$ is open in $\bigA$, so $u$ is continuous. Since $a$ was arbitrary, $a \mapsto a^{-1}$ is continuous.
\end{proof}

\begin{Dn}
The spectrum of an element $a \in \bigA$ is defined by
\[
\sigma(a)=\{\lambda \in \mathbb C : a-\lambda \text{ is not invertible}\}.
\]
\end{Dn}

In \cite[p. 167]{MR0438123} Naimark defines
topological algebras as algebras which are locally convex vector spaces, and for which the product is separately continuous
in each variable. Since strong algebras are
locally convex as strong dual complete of complete countably normed spaces,
we can apply to strong algebras which satisfy the assumption of Theorem \ref{cycontinuous}, the results of Naimark \cite[\S3, p. 169]{MR0438123} on topological algebras
with continuous inverse. In particular he showed that
\begin{Tm}[Mark Naimark]
The spectrum of any element in a locally convex unital algebra with continuous inverse is non-empty and closed.
\end{Tm}
In the next theorem, we get a bound on the spectrum, which is specific to strong algebras, even if they do not satisfy the assumption of Theorem \ref{cycontinuous}.

\begin{Tm}\label{boundspectrum}
For any $a \in \bigA$,
\[
\sigma(a) \subseteq \{z \in \mathbb C : |z|\leq \inf_{\{(p,q) : p>q+d\}}A_{p,q}\|a\|_q\}.
\]
\end{Tm}
\begin{proof}[Proof]
For every $0\neq\lambda \in \sigma(a)$,
$1- \frac{a}{\lambda}$ is not invertible and thus for any $q$ and for any $p>q+d$
\[
A_{p,q}\left\|\frac{a}{\lambda}\right\|_q \geq 1.
\]
Thus,
\[
|\lambda| \leq \inf_{\{(p,q) : p>q+d\}}A_{p,q}\|a\|_q.
\]
\end{proof}


\section{The convolution algebra $\bigcup_p L_2(S,\mu_p)$}
\setcounter{equation}{0}

\begin{Dn}
Let $G$ be a locally compact topological group with a Haar measure $\mu$ and let $S \subseteq G$ be a Borel semi-group.
The convolution of two measurable functions $f,g$ with respect to $S$ is defined by
\begin{equation}
\label{convS}
(f*g)(x)=\int_{S \cap xS^{-1}}f(y)g(y^{-1}x)d\mu(y)
\end{equation}
for any $x \in S$ such that the integral converges.
\end{Dn}

\begin{Dn}
Let $G$ be a locally compact topological group with a left Haar measure $\mu$ and let  $S \subseteq G$ be a Borel semi-group.
Let $(\mu_p)$ be a sequence of measures on $G$ such that $\mu \gg \mu_1 \gg \mu_2 \gg \cdots$.
$\bigcup_{p\in\mathbb N} L_2(S,\mu_p)$ is called a strong convolution algebra if there exists a non negative integer
$d$ such that for
every $p> q+d$
there exists a positive constant $A_{p,q}$ such that for every $f \in L_2(S,\mu_q)$ and $g \in L_2(S,\mu_p)$,
\[
\|f*g\|_p \leq A_{p,q}\|f\|_q\|g\|_p \quad \text{ and } \quad \|g*f\|_p \leq A_{p,q}\|f\|_q\|g\|_p,
\]
where $\|\cdot\|_p$ is the norm associate to $L_2(S,\mu_p)$.
\end{Dn}

\begin{Rk}
\label{rk42}
One can easily verify that requiring the measures $(\mu_p)$ of $\bigcup_p L_2(S,\mu_p)$ to obey
\[
1 \geq \frac{d\mu_1}{d\mu} \geq \frac{d\mu_2}{d\mu} \geq \cdots > 0 \quad \mu\text{-a.e.}
\]
assures that $\bigcup_p L_2(S,\mu_p)$ is a strong algebra.
\end{Rk}

The following theorem allows to define a wide family of strong convolution algebras.
For a converse theorem in the where $G$ is discrete see Theorem \ref{converse}.
\begin{Tm}\label{DFsemi}
Assume that for every $x,y \in S$ and for every $p \in \mathbb N$
\[
\frac{d\mu_p}{d\mu}(xy) \leq \frac{d\mu_p}{d\mu}(x)\frac{d\mu_p}{d\mu}(y),
\]
Then, for every choice of $f
\in L_2(S,\mu_q)$ and $g \in L_2(S,\mu_p)$,
\begin{equation}\label{ineq_Bash}
\begin{split}
\|f*g\|_{p} \leq &\left(\int_S \frac{d\mu_p}{d\mu_q}d\mu \right)^{\frac 12}\|f\|_{q}\|g\|_{p}\\
&\quad\quad\text{ and }\\
\|g*f\|_{p} \leq &\left(\int_S \frac{d\mu_p}{d\mu_q}d \widetilde \mu \right)^{\frac 12}\|f\|_{q}\|g\|_{p},
\end{split}
\end{equation}
where $\widetilde \mu$ is the right Haar measure.
In particular, if there exists a non negative integer $d$ such that
\[
\int_S \frac{d\mu_p}{d\mu_q} d\mu < \infty
\quad \text{and} \quad
\int_S \frac{d\mu_p}{d\mu_q} d\widetilde \mu < \infty
\]
for every positive integers $p,q$ such that $p > q+d$,
then $\bigcup_p L_2(S,\mu_p)$
is a strong convolution algebra (with $A_{p,q}^2=\max\left(\int_S \frac{d\mu_p}{d\mu_q} d\mu ,\int_S \frac{d\mu_p}{d\mu_q} d\widetilde\mu\right)$).
\end{Tm}

Before we prove this theorem, we need the following lemma, which is the core of the theorem.

\begin{La}\label{AG1}
Let $\nu,\lambda$ be two Borel measures on $G$ such that $\lambda
\ll \nu \ll \mu$.
If for any $x,y \in S$
\[
\frac{d\lambda}{d\mu}(xy) \leq \frac{d\lambda}{d\mu}(x)\frac{d\lambda}
{d\mu}(y)\quad \mu\text{-a.e.},
\]
then for any $f \in L_2(S,\nu)$ and $g \in L_2(S,\lambda)$
\begin{equation}
\label{ineq2}
\begin{split}
\|f*g\|_{L_2(S,\lambda)} \leq &\left(\int_S
\frac{d\lambda}{d\nu}d\mu \right)^{\frac
12}\|f\|_{L_2(S,\nu)}\|g\|_{L_2(S,\lambda)}\\
&\quad\quad \quad \text{and}\\
\|g*f\|_{L_2(S,\lambda)} \leq &\left(\int_S
\frac{d\lambda}{d\nu}d\widetilde \mu \right)^{\frac
12}\|f\|_{L_2(S,\nu)}\|g\|_{L_2(S,\lambda)}.
\end{split}
\end{equation}
where $\widetilde \mu$ is the right Haar measure.
\end{La}

\begin{proof}[Proof]
We denote
\[
f_{\lambda,\mu}(x)=f(x)\sqrt{\frac{d\lambda}{d\mu}}(x) \quad \text{ and }\quad g_{\lambda,\mu}(x)=g(x)\sqrt{\frac{d\lambda}{d\mu}}(x).
\]
Then,
\begin{align*}
\|f* g\|_{L_2(S,\lambda)} ^2
&=\int_S \left|\int_{S \cap xS^{-1}}f(y)g(y^{-1}x)d\mu(y)\right|^2d\lambda(x)\\  \displaybreak[1]
&\leq \int_S \left(\int_{S \cap xS^{-1}}|f(y)||g(y^{-1}x)|d\mu(y)\right)^2d\lambda(x)\\  \displaybreak[1]
&=\int_S \left(\int_{S \cap xS^{-1}}|f(y)||g(y^{-1}x)| \sqrt{\frac{d\lambda}{d\mu}(x)} d\mu(y)\right)^2d\mu(x)\\  \displaybreak[1]
&\leq \int_S \left(\int_{S \cap xS^{-1}}|f(y)|\sqrt{\frac{d\lambda}{d\mu}(y)}|g(y^{-1}x)| \sqrt{\frac{d\lambda}{d\mu}(y^{-1}x)} d\mu(y)\right)^2d\mu(x)\\  \displaybreak[1]
&= \int_S \left(\int_{S \cap xS^{-1}} |f_{\lambda,\mu}(y)||g_{\lambda,\mu}(y^{-1}x)| d\mu(y)\right)^2d\mu(x)\\  \displaybreak[1]
&=\int_S \int_{S \cap xS^{-1}} \int_{S \cap xS^{-1}}|f_{\lambda,\mu}(y)||f_{\lambda,\mu}(\tilde y)|| g_{\lambda,\mu}(y^{-1}x)||g_{\lambda,\mu}(\tilde y^{-1}x)| d\mu(y)d\mu(\tilde y) d\mu(x)\\  \displaybreak[1]
&=\int_S \int_S| f_{\lambda,\mu}(y)||f_{\lambda,\mu}(\tilde y)|\left( \int_{S \cap yS \cap \tilde yS} |g_{\lambda,\mu}(y^{-1}x) ||g_{\lambda,\mu}(\tilde y^{-1}x)|d\mu(x)\right)d\mu(y)d\mu(\tilde y) \\ \displaybreak[1]
&\leq \int_S \int_S |f_{\lambda,\mu}(y)||f_{\lambda,\mu}(\tilde y)| \left\|g_{\lambda,\mu}\right\|_{L_2(S,\mu)}^2 d\mu(y)d\mu(\tilde y)\\  \displaybreak[1]
&=\left(\int_S |f_{\lambda,\mu}(y)| d\mu(y) \right)^2 \left\|g\right\|_{L_2(S,\lambda)}^2\\ \displaybreak[1]
&=\left|\int_S  \sqrt{\frac{d\lambda}{d\nu}}| f_{\lambda,\nu}| d\mu \right|^2 \left\|g\right\|_{L_2(S,\lambda)}^2\\  \displaybreak[1]
&\leq \left(\int_S \frac{d\lambda}{d\nu}d\mu \right) \left\|f_{\nu,\mu}(y)\right\|_{L_2(S,\mu)}^2\|g\|_{L_2(S,\lambda)}^2 \\ \displaybreak[1]
&=\left(\int_S \frac{d\lambda}{d\nu}d\mu \right) \left\|f \right\|_{L_2(S,\nu)}^2\|g\|_{L_2(S,\lambda)}^2
\end{align*}
yields the first inequality of \eqref{ineq2}. As for the second inequality, note that
\[
(g*f)(x)=\int_{S \cap xS^{-1}}g(y)f(y^{-1}x)d\mu(y)=\int_{S \cap S^{-1}x}f(y)g(xy^{-1})d\widetilde{\mu}(y).
\]
Thus, replacing the terms $\mu$, $S \cap xS^{-1}$ and $y^{-1}x$ in the proof of the first inequality, with the terms
$\tilde \mu$, $S \cap S^{-1}x$ and $xy^{-1}$ respectively, yields a proof for the second inequality. .
\end{proof}

We are now ready to prove Theorem \ref{DFsemi}.
\begin{proof}[Proof of Theorem \ref{DFsemi}]
In view of Lemma \ref{AG1}, it holds that
\[
\begin{split}
\|f*g\|_{p} \leq &\left(\int_S \frac{d\mu_p}{d\mu_q}d\mu \right)^{\frac 12}\|f\|_{q}\|g\|_{p}\\
&\quad\quad\text{ and }\\
\|g*f\|_{p} \leq &\left(\int_S \frac{d\mu_p}{d\mu_q}d \widetilde \mu \right)^{\frac 12}\|f\|_{q}\|g\|_{p},
\end{split}
\]
for any $p>q+d$. This yields the requested result.
\end{proof}

In the following theorems we focus on the case where $G$ is a discrete group. In this case, the Haar measure is simply the counting measure. Nonetheless, we keep using the notations of the general case, namely integrals instead of sums, and $L_2$ instead of $\ell_2$.
\begin{Tm}\label{converse}
Let $\bigcup_p L_2(S,\mu_p)$ is a strong convolution algebra, where $S$ is a semi-group in a discrete group $G$, such that
\[
\frac{d\mu_p}{d\mu}=\left(\frac{d\mu_1}{d\mu}\right)^p \text { for every } p \in \mathbb N,
\]
and such that there exists $d$ for which $\int_S\left( \frac{d\mu_1}{d\mu}\right)^{d}d\mu < \infty$. Then for any $x,y \in S$,
\[
\frac{d\mu_p}{d\mu}(xy) \leq \frac{d\mu_p}{d\mu}(x)\frac{d\mu_p}{d\mu}(y) \text { for every } p \in \mathbb N.
\]
\end{Tm}
\begin{proof}
Denoting by $\delta_x$ the characteristic function of $\{x\}$, we have
\[
\|\delta_x*\delta_y\|_{p} \leq \left(\int_S\left( \frac{d\mu_1}{d\mu}\right)^{p-q}d\mu \right)^{\frac 12}\|\delta_x\|_{q}\|\delta_y\|_{p}.
\]
Setting $q=p-(d+1)$ we obtain
\[
\frac{d\mu_p}{d\mu}(xy) \leq \left(\int_S\left( \frac{d\mu_1}{d\mu}\right)^{d+1} d\mu\right)^{\frac 12}\frac{d\mu_{p-(d+1)}}{d\mu}(x)\frac{d\mu_p}{d\mu}(y),
\]
or
\[
\frac{d\mu_1}{d\mu}(xy) \leq \left(\int_S\left( \frac{d\mu_1}{d\mu}\right)^{d+1} d\mu\right)^{\frac 1{2p}} \left(\frac{d\mu_{1}}{d\mu}(x)\right)^{\frac{p-(d+1)}{p}}\left(\frac{d\mu_1}{d\mu}(y)\right).
\]
Taking $p \to \infty$ yields the requested result.
\end{proof}

\begin{Tm}\label{nuclear}
If $G$ is discrete and for $p>q$, $\int_{S}\frac{d\mu_p}{d\mu_q}d\mu<\infty$, then the embedding
\[
T_{q,p}:L_2(S,\mu_q) \hookrightarrow L_2(S,\mu_q)
\]
is Hilbert-Schmidt. As a result, if for any $p>q+d$ the above integral converges (as in the assumption of Theorem \ref{DFsemi}), then $\bigcup_p L_2(S,\mu_p)$ is nuclear.
\end{Tm}
\begin{proof}[Proof]
For any $x\in S$, $\|\delta_x\|_q^2=\frac{d\mu_q}{d\mu}(x)$. Thus, $e_x^{(q)}=\left(\frac{d\mu_q}{d\mu}(x)\right)^{-\frac 12}\delta_x$ ($x \in S$) form an orthonormal basis of $L_2(S,\mu_q)$, and
\[
\|T_{q,p}\|_{HS}^2=\sum_{x\in S}\|T_{q,p}e_x^{(q)}\|_p^2=\int_S\frac{d\mu_p}{d\mu_q}d\mu,
\]
where $\|\cdot\|_{HS}$ denotes the Hilbert Schmidt norm of the embedding $T_{q,p}$.
\end{proof}

A first example of a strong convolution algebra was given in our previous paper \cite{vage1}. We briefly discuss it now.
Recall first that
the Schwartz space of complex tempered distributions can be viewed as the union of the Hilbert spaces
\begin{equation*}
s'_p:=\left\{a \in \mathbb C^{\mathbb N_0}: \sum_{n=0}^\infty (n+1)^{-2p}|a(n)|^2< \infty \right\}.
\end{equation*}
Let $G=\mathbb Z$ with its discrete topology, $\mu$ the counting measure (which is also the Haar measure of $G$), $S=\mathbb N_0$ and setting
\begin{equation}
\label{sprime}
d\mu_p(n)=(n+1)^{-2p},
\end{equation}
we conclude that $s'_p=L_2(S,\mu_p)$, and $s'=\bigcup_{p \in \mathbb N}s'_p=\bigcup_{p \in \mathbb N}L_2(S,\mu_p)$.
The convolution then becomes the standard one sided convolution of sequences.\\

One may ask whether $s'$ is a strong convolution algebra, i.e. if there exists a constant $d$ such that for
any $p> q+d$ there exists a positive constant $A_{p,q}$ such that for any $a \in s'_q$ and $b \in s'_p$,
\[
\|a*b\|_{p} \leq A_{p,q}\|a\|_{q}\|b\|_{p} \quad \text{ and }\quad  \|b*a\|_{p} \leq A_{p,q}\|a\|_{q}\|b\|_{p}.
\]
Since,
\[
\int_S\left( \frac{d\mu_1}{d\mu}\right)^{d}d\mu=\sum_{n=0}^\infty (n+1)^{-2} < \infty
\]
and
\[
\frac{d\mu_p}{d\mu}(n)=(n+1)^{-2p}=\left(\frac{d\mu_1}{d\mu}\right)^p \text { for any } p \in \mathbb N,
\]
we conclude in view of the last theorem that the answer is negative, that is, $s'=\bigcup_{p \in \mathbb N}s'_p$ is not a
strong convolution algebra. In \cite{vage1} we replace the measures $(n+1)^{-2p}$ by $2^{-np}$, and
obtain a strong convolution algebra that contains $s'$, and which can be identified as the dual of a space of
entire holomorphic functions that is included in the Schwartz space of complex-valued rapidly decreasing
smooth functions.\\

Other examples of strong convolution algebras, which can be constructed in the manner described in the last theorem,
are given in Sections \ref{sec5} and \ref{sec6} respectively.\\

To end this section, we make a brief discussion on the case where
$G$ is non-unimodular (i.e. its left and right Haar measure do
not coincide). As is clear from the statement of Theorem
\ref{DFsemi} and from the proof of Lemma \ref{AG1}, if $G$ is a
non-unimodular group, since the integrals $\int_S
\frac{d\mu_p}{d\mu_q} d\mu ,\int_S \frac{d\mu_p}{d\mu_q}
d\widetilde\mu$ need not be equal, it may happen that one of the
inequalities \eqref{ineq_Bash} is stricter than the other.
Nonetheless, if both integrals are finite, then we may always
take $A_{p,q}^2=\max\left(\int_S \frac{d\mu_p}{d\mu_q} d\mu
,\int_S \frac{d\mu_p}{d\mu_q} d\widetilde\mu\right)$, and obtain
\[
\|f*g\|_p \leq A_{p,q}\|f\|_q\|g\|_p \quad \text{ and } \quad \|g*f\|_p \leq A_{p,q}\|f\|_q\|g\|_p.
\]
As an example, consider the non-unimodular group $G=\{(a,b)\in \mathbb R^2 \,:\,a>0\}$ with the multiplication
$(a,b)\cdot(c,d)=(ac,ad+b)$ which can be identified as the subgroup of $GL_2(\mathbb R)$
\[
\left\{
\begin{pmatrix}
a&b\\
0&1
\end{pmatrix}
\in GL_2(\mathbb R):a>0\right\}.
\]
This is the so-called $ax+b$ group. One can easily verify that its left Haar measure is given by $d\mu(a,b)=a^{-1}dadb$, and that its right Haar measure is given by $d\widetilde{\mu}(a,b)=a^{-2}dadb$.
Let $S$ be the semigroup
\[
\{(a,b)\in \mathbb R^2 \,:\,a \geq 1,b\geq 0 \}
\]
and set $d\mu_p(a,b)=a^{-(p+2)}e^{-bp}dadb$. So $\bigcup_p L_2(S,\mu_p)$ is the space of all measurable functions $S \to \mathbb C$ satisfy
\[
\int_0^\infty \int_1^\infty |f(a,b)|^2a^{-(p+2)}e^{-bp}dadb < \infty \quad \text{for some }p\in \mathbb N,
\]
with a convolution product
\[
(f*g)(a,b)=\int_0^b \int_1^a f(x,y)g \left(x^{-1}a,x^{-1}(b-y)\right)a^{-1}dadb.
\]
Since clearly $\frac{d\mu_p}{d\mu}\left((a,b)(a',b')\right)\leq \frac{d\mu_p}{d\mu}(a,b)\frac{d\mu_p}{d\mu}(a',b')$,
and since for any $p>q$,
\[
\int_S \frac{d\mu_p}{d\mu_q} d\mu=\frac{1}{(2(p-q)+1)(p-q)}\quad \text{ and }\quad \int_S \frac{d\mu_p}{d\mu_q} d\widetilde\mu=\frac{1}{2(p-q)^2},
\]
By theorem \ref{DFsemi}, $\bigcup_p L_2(S,\mu_p)$ is a strong convolution algebra with
$
A_{p,q}=\frac 1 {\sqrt{2}(p-q)}.
$
However, in view of \eqref{ineq_Bash}, on of the inequalities is stricter. More precisely, denoting $B_{p,q}=\left(\frac{1}{(2(p-q)+1)(p-q)}\right)^{\frac 12}$,
for every $f \in L_2(S,\mu_q)$ and $g \in L_2(S,\mu_p)$ ($p>q$) we have,
\[
\|f*g\|_p \leq B_{p,q}\|f\|_q\|g\|_p \quad \text{ and } \quad \|g*f\|_p \leq A_{p,q}\|f\|_q\|g\|_p.
\]
A question which we leave open is whether there exists a ``one-sided strong convolution algebra" which is not a (``two-sided") strong convolution algebra, i.e. an algebra of the form
$\bigcup_{p\in\mathbb N} L_2(S,\mu_p)$ of which there are $d$ and $A_{p,q}$ such that
 for every $p>q+d$ and for every $f \in L_2(S,\mu_q)$ and $g \in L_2(S,\mu_p)$,
$
\|f*g\|_p \leq A_{p,q}\|f\|_q\|g\|_p,
$
but where the ``reflected" inequality,
namely $\|g*f\|_p \leq A_{p,q}\|f\|_q\|g\|_p$,
does not hold for any choice of $d$ and $A_{p,q}$.
Clearly, such an example should be over a non-unimodular group.\\


\section{The space of germs of  holomorphic functions in zero}
\setcounter{equation}{0}
Let $H(\mathbb C)$ be the space of entire holomorphic functions. It can be topologized as follows.
For any $p \in \mathbb N$ we define the $n$-norm on $H(\mathbb C)$ by
\[
\|f\|_p^2=\frac{1}{2\pi}\int_{|z|=2^p}\left|f(z)\right|^2dz,
\]
then the topology of $H(\mathbb C)$ is the associated Fr\'echet topology.
\begin{Pn}
The topology defined above coincides with the usual topology of $H(\mathbb C)$, that is the topology of uniform convergence on compact sets.
\end{Pn}
\begin{proof}[Proof]
Clearly,
\[
\|f\|_p^2=\frac{1}{2\pi}\int_{|z|=2^p}\left|f(z)\right|^2dz \leq 2^p \sup_{z \in 2^p \Bar{\mathbb D}}|f(z)|^2.
\]
On the opposite direction, using Cauchy theorem, for any $z \in 2^p \Bar{\mathbb D}$
\[
\begin{split}
\left|f(z)\right|^2
&=\left|\frac{1}{2\pi}\int_{|z|=2^{p+1}}\frac{f(\omega)^2}{\omega-z}d\omega\right| \\
&\leq \frac{1}{2\pi}\int_{|z|=2^{p+1}}\left|\frac{f(\omega)^2}{\omega-z}\right|d\omega \\
&\leq 2^{-p}\frac{1}{2\pi}\int_{|z|=2^{p+1}}\left|f(\omega)\right|^2d\omega \\
&=2^{-p}\|f\|_{p+1}^2
\end{split}
\]
Thus,
\[
\sup_{z \in 2^p \bar{\mathbb D}}|f(z)| \leq 2^{-p/2}\|f\|_{p+1},
\]
and we conclude that
\[
\|f\|_p \leq 2^{p/2} \sup_{z \in 2^p \bar{\mathbb D}}|f(z)| \leq \|f\|_{p+1}.
\]
\end{proof}

Let ${\bf H_2}(2^p \mathbb D)$ be the Hardy space of the disk $2^p \mathbb D$, i.e.
\[
{\bf H_2}(2^p \mathbb D)=\left\{ f \text{ is holomorphic on } 2^p \mathbb D : \sup_{0<r<1}\frac{1}{2\pi}\int_{|z|=r\cdot 2^p}\left|f(z)\right|^2dz < \infty \right\},
\]
then clearly, ${\bf H_2}(2^p \mathbb D)$ is the completion of $H(\mathbb C)$ with respect to the $p$-norm, and
\[
{\bf H_2}(\mathbb D) \supseteq {\bf H_2}(2 \mathbb D) \supseteq \cdots \supseteq {\bf H_2}(2^p \mathbb D) \supseteq \cdots.
\]

Therefore, we obtain
\[
H(\mathbb C)=\bigcap_{p \in \mathbb N} {\bf H_2}(2^p \mathbb D),
\]
i.e. $H(\mathbb C)$ is a countably Hilbert space.\\

The Hardy space can be viewed also in terms of power series, and with
$f(z)=\sum_{n=0}^\infty a_n z^n$ , it holds that
\[
{\bf H_2}(2^p \mathbb D)=\left\{ \sum_{n=0}^\infty a_n z^n: \sum_{n=0}^\infty |a_n|^2 2^{2np}< \infty \right\},
\]
and
\[
 \sup_{0<r<1}\frac{1}{2\pi}\int_{|z|=r\cdot 2^p}\left|f(z)\right|^2dz  =  \sum_{n=0}^\infty |a_n|^2 2^{2np}.
\]
Hence we obtain,
\[
H(\mathbb C)=\left\{ \sum_{n=0}^\infty a_n z^n: \sum_{n=0}^\infty |a_n|^2 2^{2np}< \infty \text { for all } p \in \mathbb N\right\}
\]
with the associated topology.\\

Now, let $H'(\mathbb C)$ denote the dual space of $H(\mathbb C)$,
equipped with the strong topology as dual of countably Hilbert
space. Then,
\[
H'(\mathbb C)=\bigcup_{p \in \mathbb N} {\bf H_2}(2^{-p} \mathbb D),
\]
where it holds that
\[
{\bf H_2}(\mathbb D) \subseteq {\bf H_2}(2^{-1} \mathbb D) \subseteq \cdots \subseteq {\bf H_2}(2^{-p} \mathbb D) \subseteq \cdots.
\]
Therefore, we conclude
\begin{Pn}
$H'(\mathbb C)$ can be identified with the space of holomorphic germs at zero $\mathcal O_0$.
\end{Pn}

The multiplication of two elements in $H'(\mathbb C)$,
$f(z)=\sum_{n=0}^\infty a_nz^n$ and $g(z)=\sum_{n=0}^\infty
b_nz^n$ is defined by
\[
(fg)(z)=\sum_{n=0}^\infty \left(\sum_{m=0}^n a_m b_{n-m} \right) z^n.
\]
\begin{Tm}\label{SAgerms}
Let $f \in {\bf H_2}(2^{-q} \mathbb D)$ and $g \in {\bf H_2}(2^{-p} \mathbb D)$ where $p>q+d$.
Then $fg \in {\bf H_2}(2^{-p} \mathbb D)$. Furthermore,
\[
\|fg\|_{-p} \leq (1-4^{-(p-q)})^{-1}\|f\|_{-q}\|g\|_{-p} .
\]
\end{Tm}
\begin{proof}[Proof]
This follows directly from Theorem \ref{DFsemi}.
\end{proof}
In view of Theorems \ref{SAcont}, \ref{cycontinuous},
\ref{boundspectrum} and \ref{nuclear}, we conclude:
\begin{Cy}\label{TopGerms}
The following properties hold:
\begin{enumerate}[(a)]
\item The multiplication is continuous, i.e. $(H'(\mathbb C),+,\cdot)$ is a topological algebra.
\item $GL(H'(\mathbb C))$ is open and the function $f \mapsto f^{-1}$ is continuous.
\item $f \in H'(\mathbb C)$ is invertible if and only if $f(0) \neq 0$.
\item $H'(\mathbb C)$ is nuclear.
\end{enumerate}
\end{Cy}

Note that the last item is well known; see for instance \cite[pp.
105-106]{pietsch} (since the dual of nuclear space is nuclear) or
\cite{gr1,gr2}.

\begin{Rk}
Let $(\xi_n)$ denote the Hermite functions, which form an orthonomal basis of $L_2(\mathbb C)$, and
consider the isometrically isomorphism $\mathcal G:{L}^2(\mathbb R) \to {\bf H_2}(\mathbb D)$ defined by
\[
\mathcal G: \xi_n \mapsto z^n.
\]
In \cite{vage1} we showed that the space $\mathscr G$ of all entire functions $f(z)$ such that
\[
 \iint_\mathbb{C} \left|f(z)\right|^2 e^{\frac{1-2^{-2p}}{1+2^{-2p}}x^2
 -\frac{1+2^{-2p}}{1-2^{-2p}}y^2}dxdy<\infty \quad \text { for all } p \in \mathbb N.
\]
can be viewed as the space of all functions $\sum_{n=0}^\infty a_n \xi_n$ such that for any $p \in \mathbb N$
\[
\sum_{n=0}^\infty |a_n|^2 2^{2np} <\infty.
\]
Thus, the image of $\mathscr G$ under $\mathcal G$ is $H'(\mathbb C)$.
We also note that $\mathscr G$ is included in the Schwartz space of rapidly decreasing smooth complex functions $\mathscr S$. Thus, the Schwartz space of tempered distributions $\mathscr S'$ is included in  $\mathscr G'$ which can be identified with the space of gems of holomorphic functions in zero.
\end{Rk}


\section{The space $\bigP$}\label{sec5}
\setcounter{equation}{0}
We now present a new example of a convolution algebra $\bigcup_p L_2(S,\mu_p)$.
Let $\bigP$ be the space of all functions $f:[1,\infty) \to \mathbb R$ such that there exists $p \in \mathbb N$  with
\[
\int_1^\infty |f(x)|^2 \frac{dx}{x^{p+1}}<\infty.
\]
In particular, any restriction of a polynomial function into $[1,\infty)$ belongs to $\bigP$. Thinking of $[1,\infty)$
as a group with respect to the multiplication of the real numbers, and since the Haar measure
of $((0,\infty), \cdot)$ is $\frac{dx}{x}=d(\ln(x))$, we obtain the following convolution
\[
(f*g)(x)=\int_1^xf(y)g\left(\frac{x}{y}\right)\frac{dy}{y}, \quad \forall x \in [1,\infty).
\]
We also note that for any $p>q$,
\[
\int_1^\infty\frac{x^{-{(p+1)}}}{x^{-{(q+1)}}} \frac{dx}{x}=\frac{1}{p-q}.
\]
Thus, using Theorem \ref{DFsemi}, we obtain that for any $p>q$,
\begin{equation}\label{P_ineq}
\int_1^\infty \left( \int_1^x f(y)g\left(\frac{x}{y}\right)\frac{dy}{y}\right)^2 \frac{dx}{x^{p+1}}
\leq \frac{1}{p-q} \int_1^\infty f^2(x)\frac{dx}{x^{q+1}}\int_1^\infty g^2(x)\frac{dx}{x^{p+1}}.
\end{equation}

We now present inequalities on Dirichlet series. To that purpose we note that in
\eqref{P_ineq} we can assume that $p$ and $q$ are positive real numbers.
\begin{Dn}
The one-sided Mellin transform is defined by
\[
(\bigM f)(s)=\int_1^\infty x^s f(x) \frac{dx}{x}.
\]
\end{Dn}
Thus, \eqref{P_ineq} can be rewritten as
\begin{equation}\label{P_ineq2}
(\bigM (f*g)^2)(-t) \leq \frac{1}{t-r} (\bigM f^2)(-r)(\bigM g^2)(-t) \quad \text{for any }t>r.
\end{equation}
\begin{Ex}
Consider the Dirichlet series
\[
\sum_{n=1}^\infty a_n n^{-s} \quad \text{and } \quad \sum_{n=1}^\infty b_n n^{-s}.
\]
Since for any $s$ in the half-plane of absolute convergence,
\[
\sum_{n=1}^\infty a_n n^{-s} =s \int_1^\infty y^{-s} \sum_{n \leq y}a_n \frac{dy}{y}=s\left(\bigM \sum_{n \leq y}a_n\right)(-s),
\]
(see \cite[Theorem 1.3, p. 13]{MR2378655}, \cite[p. 10]{MR1433800})
\eqref{P_ineq2} yields
\[
\int_1^\infty \left( \int_1^x \sqrt{\sum_{n \leq y}a_n \sum_{n \leq \frac{x}{y}}b_n}\frac{dy}{y}\right)^2 \frac{dx}{x^{t+1}}
\leq \frac{1}{t(t-r)r} \sum_{n=1}^\infty a_n n^{-r} \sum_{n=1}^\infty b_n n^{-t}.
\]
For example taking the zeta function of Riemann,
\[
\zeta(s)=\sum_{n=1}^\infty n^{-s},
\]
which converges for ${\rm Re}~s>1$, one obtains
\[
\int_1^\infty \left( \int_1^x \sqrt{\ \lfloor y \rfloor\left \lfloor\frac{x}{y} \right \rfloor}\frac{dy}{y}\right)^2 \frac{dx}{x^{t+1}}
\leq \frac{\zeta(t) \zeta(r)}{t(t-r)r} \quad \text{ for any } 1<r<t .
\]
Hence,
\[
\int_1^\infty \left( \int_1^x \sqrt{\ \lfloor y \rfloor\left \lfloor\frac{x}{y} \right \rfloor}\frac{dy}{y}\right)^2 \frac{dx}{x^{t+1}}
\leq \inf_{r \in (1,t)}\left(\frac{\zeta(t) \zeta(r)}{t(t-r)r} \right)\quad \text{ for any } t>1 .
\]
In a similar way, if we take $\varphi(n)$ to be the phi Euler function (that is $\varphi(n)$ is the number of positive
integers less than or equal to $n$ and relatively prime with $n$), we obtain
\[
\int_1^\infty \left( \int_1^x \sqrt{\sum_{n \leq y}\varphi(n) \sum_{n \leq \frac{x}{y}}\varphi(n)}\frac{dy}{y}\right)^2 \frac{dx}{x^{t+1}}
\leq \inf_{r \in (1,t)}\left(\frac{\zeta(t-1) \zeta(r-1)}{t\zeta(t)(t-r)r\zeta(r)} \right).
\]
\end{Ex}


\section{The Kondratiev space of stochastic distributions}
\label{sec6}
\setcounter{equation}{0}
We now present two other examples of strong convolution algebras. First we consider
the Kondratiev space ${\mathcal S}_{-1}$ of stochastic
distributions. The space ${\mathcal S}_{-1}$
plays an important role in white noise space analysis and in the theory of linear stochastic differential equations and
linear stochastic systems; see \cite{MR1408433,alp,MR2610579,vage1, aa_goh}.
Next we introduce its non commutative counterpart $\tilde {\mathcal S}_{-1}$. Further properties of $\tilde
{\mathcal S}_{-1}$ and applications will be presented in a forthcoming publication.\\

Let
\begin{equation}
\ell= \mathbb N_0^{(\mathbb N)}=\left\{\alpha \in \mathbb N_0^{\mathbb N} : \text{supp}(\alpha) \text{ is finite} \right\}=\oplus_{n\in \mathbb N}\mathbb N_0 e_n,
\end{equation}
and let $\bigS_{-1}$ be the space of all functions $f:\ell \to \mathbb C$ such that there exists $p \in \mathbb N$  with
\[
\sum_{\alpha\in \ell} |f(\alpha)|^2(2\mathbb N)^{-\alpha p}<\infty.
\]
Denoting
\[
\mu_p(\alpha)=(2\mathbb N)^{-\alpha p}=\prod_{n=1}^\infty (2n)^{-\alpha(n) p}
\]
we conclude that $\bigS_{-1}$ is the convolution algebra $\bigcup_p L_2(S,\mu_p)$,
with the convolution
\[
(f*g)(\alpha)=\sum_{\beta \leq \alpha} f(\beta)g(\alpha-\beta), \quad \forall \alpha \in \ell.
\]
In stochastic analysis, this convolution is called the Wick product; see \cite{MR2387368}.
We also note that for every $p> q+1$,
\[
\sum_{\alpha \in \ell} \frac{(2\mathbb N)^{-\alpha p}}{(2\mathbb N)^{-\alpha q}}=\sum_{\alpha \in \ell}
(2\mathbb N)^{-\alpha (p-q)}=A(p-q)^2<\infty.
\]
The fact that the last term is finite is due to \cite{zhang}. Another proof can be found in \cite{vage1}.
Thus, due to Theorem \ref{DFsemi}, we obtain that for every $p>q+1$,
\begin{equation}\label{S1_ineq}
\|f*g\|_p \leq A(p-q)\|f\|_q\|g\|_p.
\end{equation}
This result (in the context of Kondratiev space of stochastic distributions) was first proved by V\aa ge in 1996
(see \cite{vage96}).\\

We now introduce a construction of the non-commutative version of $\bigS_{-1}$ as a convolution algebra
(in a forthcoming paper we show that this space can also be obtained via a union of full Fock spaces).
We replace the free commutative semi-group generated by $\mathbb N$,
namely $\ell$, with the free (non-commutative) semi-group generated by $\mathbb N$, which we will denote by $\tilde \ell$.
For $\alpha=z_{i_1}^{\alpha_{1}}z_{i_2}^{\alpha_{2}} \cdots z_{i_n}^{\alpha_{n}} \in \tilde \ell$
(where $i_1 \neq i_2 \neq \cdots \neq i_n$) we define
\[
(2 \mathbb N)^\alpha=\prod_{k=1}^n(2i_k)^{\alpha_k}=\prod_{j \in \{i_1,\dots,i_n\}}(2j)^{\left( \sum_{
k\in\{m\, :\, i_m=j\}}\alpha_k \right)}.
\]

\begin{Pn} For every $p,q$ integers such that $p> q+1$, it holds that
\[
B(p-q)^2=\sum_{\alpha \in \tilde \ell} (2 \mathbb N)^{-\alpha(p-q)}<\infty.
\]
\end{Pn}
\begin{proof}[Proof]
We note that
\[
\begin{split}
B(p-q)^2
&=\sum_{\alpha \in \tilde \ell} (2 \mathbb N)^{-\alpha(p-q)}\\
&=\sum_{n=0}^\infty\sum_{\alpha \in \tilde \ell,|\alpha|=n}\prod_{k=1}^\infty(2k)^{-\alpha_k(p-q)}\\
&=\sum_{n=0}^\infty\sum_{\alpha \in \ell,|\alpha|=n}\frac{n!}{\alpha!}\prod_{k=1}^\infty(2k)^{-\alpha_k(p-q)}.
\end{split}
\]
Considering now an experiment with $\mathbb N$ results, where the probability of the result $k$ is
$p_k=c \cdot (2k)^{-(p-q)}$ ($c$ is chosen such that $\sum p_k=1$, and clearly for any $p > q+1$, $c>1$),
the probability that repeating the experiment $n$ times yields that the result $k$ occurs $\alpha_k$ times for any $k$ is
\[
\frac{n!}{\alpha!}\prod_{k=1}^\infty p_k^{\alpha_k}
=c^{n}\frac{n!}{\alpha!}\prod_{k=1}^\infty (2k)^{-\alpha_k(p-q)}.
\]
Thus,
\[
\sum_{\alpha \in \ell,|\alpha|=n}\frac{n!}{\alpha!}\prod_{k=1}^\infty(2k)^{-\alpha_k(p-q)}=c^{-n},
\]
and we obtain the requested result.
\end{proof}
We conclude that the non-commutative version of the Kondratiev space is also a strong convolution algebra.


\bibliographystyle{plain}
\def\cprime{$'$} \def\lfhook#1{\setbox0=\hbox{#1}{\ooalign{\hidewidth
  \lower1.5ex\hbox{'}\hidewidth\crcr\unhbox0}}} \def\cprime{$'$}
  \def\cfgrv#1{\ifmmode\setbox7\hbox{$\accent"5E#1$}\else
  \setbox7\hbox{\accent"5E#1}\penalty 10000\relax\fi\raise 1\ht7
  \hbox{\lower1.05ex\hbox to 1\wd7{\hss\accent"12\hss}}\penalty 10000
  \hskip-1\wd7\penalty 10000\box7} \def\cprime{$'$} \def\cprime{$'$}
  \def\cprime{$'$} \def\cprime{$'$}

\end{document}